\documentclass[12pt]{amsart}

\usepackage{amsmath}
\usepackage[mathscr]{eucal}
\usepackage{amssymb}
\usepackage{latexsym}
\usepackage{amsthm}
\usepackage{amscd}
\usepackage{txfonts}

\theoremstyle{plain}

\newtheorem{theorem}{Theorem}[section]

\newtheorem{lemma}[theorem]{Lemma}

\newtheorem*{main}{{\bf Theorem}}

\theoremstyle{definition}

\newtheorem{remark}[theorem]{Remark}

\newcommand{\reflem}[1]{Lemma~{\rm \ref{#1}}}

\newcommand{\reftable}[1]{Table~{\rm \ref{#1}}}

\newcommand{\ol}[1]{\overline{#1}}

\newcommand{\bQ}{\mathbb{Q}}
\newcommand{\bZ}{\mathbb{Z}}
\newcommand{\bF}{\mathbb{F}}
\newcommand{\bP}{\mathbb{P}}
\newcommand{\bC}{\mathbb{C}}

\newcommand{\bFp}{\mathbb{F}_p}

\newcommand{\bFq}{\mathbb{F}_q}

\newcommand{\cL}{\mathcal{L}}
\newcommand{\cK}{\mathcal{K}}

\newcommand{\SL}{\mathrm{SL}}

\newcommand{\Rank}{\mathrm{rank}\,}

\newcommand{\Char}{\mathrm{Char}\,}

\newcommand{\Cf}{\textrm{cf.}\;}

\newcommand{\ra}{\rightarrow}

\newcommand{\lan}{\langle}
\newcommand{\ran}{\rangle}
\newcommand{\lv}{\lvert}
\newcommand{\rv}{\rvert}

\newcommand{\fg}{\pi_1}
\newcommand{\isom}{\tilde{\ra}}

\newcommand{\gl}{\lambda}
\newcommand{\gm}{\mu}

\newcommand{\disp}[1]{\displaystyle{#1}}

\newcommand{\vphi}{\varphi}

\newcommand{\bFpn}{\mathbb{F}_{p^n}}

\newcommand{\ga}{\alpha}

\newcommand{\ssm}{\smallsetminus}

\newcommand{\quadres}[2]{\begin{pmatrix} \dfrac{#1}{#2} \end{pmatrix}}
\newcommand{\dquadres}[2]{\begin{pmatrix} \dfrac{#1}{#2} \end{pmatrix}}
\newcommand{\tquadres}[2]{\begin{pmatrix} \tfrac{#1}{#2} \end{pmatrix}}

\title{Hasse-Weil zeta functions of $\SL_2$-character varieties of arithmetic two-bridge link complements}
\author{Shinya Harada}

\keywords{ 
 Hasse-Weil zeta function, $\SL_2$-character variety,
 arithmetic two-bridge link}

\begin{document}

\pagestyle{myheadings}
\pagenumbering{arabic}

\maketitle

\markboth{Shinya Harada}{Zeta of Canonical components of arithmetic two-bridge link groups}

\begin{abstract}
 Hasse-Weil zeta functions of
 $\SL_2$-character varieties of arithmetic two-bridge link groups are determined.
 Special values of the zeta functions at $s=0,1,2$
 are also investigated.
\end{abstract}

\section{Introduction}


%
%

 In \cite{HWfigHarada}
 the Hasse-Weil zeta function of
 the $\SL_2$-character variety of the figure $8$ knot,
 the only arithmetic knot in the $3$-sphere,
 has been determined
 (see \cite{BookMaRe} as a reference on arithmetic manifolds).
 In this note we extend the result to
 the most basic class of arithmetic
 $2$-component links,
 the family of arithmetic two-bridge links,
 $L_0:=5^2_1=(8/3)$ (Whitehead link),
 $L_1:=6^2_2=(10/3)$ and $L_2:=6^2_3=(12/5)$.
 In fact,
 for the canonical components (irreducible components
 containing holonomy representations) of
 their $\SL_2(\bC)$-character varieties
 which are algebraic surfaces defined by
\begin{align*}
 f_0&:=z^3 -xyz^2+(x^2+y^2-2)z -xy,\\
 f_1&:=z^4-xyz^3+(x^2+y^2-3)z^2-xyz+1,\\
 f_2&:=z^3-xyz^2+(x^2+y^2-1)z-xy,
\end{align*}
%
%
 we compute the Hasse-Weil zeta functions
$$
\zeta(L_i,s):=\zeta(f_i,s):=\prod_p \exp\left(\sum_{n=1}^{\infty}\dfrac{\#V(f_i)(\bFpn)}{n} (p^{-s})^n\right),
$$
 where $p$ runs through all the prime numbers
 and
 we denote by $V(f_i)(\bFpn)$ the set of $\bFpn$-rational
 points of $f_i$:
$$
 V(f_i)(\bFpn):= \left\{(a,b,c) \in (\bFpn)^3\mid f_i(a,b,c) = 0 \right\}
$$
 for the finite field $\bFpn$ with $p^n$ elements.
 The result is as follows.

\begin{main}
 The Hasse-Weil zeta functions
 of the canonical components of $L_0=5^2_1=(8/3)$, $L_1=6^2_2=(10/3)$ and $L_2=6^2_3=(12/5)$
 are the following.
\begin{align*}
 \zeta(L_0,s)&
%
%
%
 =\zeta_{\bQ(\sqrt{2})}(s-1) \zeta(s)^2 \zeta(s-1)^2 \zeta(s-2) (1 - 2^{1-s})^3 (1 - 2^{-s}).\\
%
\zeta(L_1,s)&
%
             =
 \zeta_{\bQ(\sqrt{5})}(s-1)^2 \zeta(s)^3 \zeta(s-2)
(1 - 2^{1-s})^{3}(1 + 2^{1-s})(1 - 2^{-s}).\\
%
 \zeta(L_2,s)&
             = \zeta(s)^2 \zeta(s-2) (1-2^{-s}).
\end{align*}
\end{main}
\noindent
 Here
 $\zeta_K(s)$ is the Dedekind zeta function
 of an algebraic number field $K$
 and $\zeta(s):=\zeta_{\bQ}(s)$ is the Riemann zeta function.
%
%
 It is shown in \cite{EmilyLandes,CanATBL} that
 the canonical components of the $\SL_2(\bC)$-character
 varieties of arithmetic two-bridge links
 have the structure of conic bundles over
 the projective line $\bP^1$.
 That makes it easy to compute the number of rational
 points of the canonical components over finite fields.
 Note that
 the fields $\bQ(\sqrt{2})$, $\bQ(\sqrt{5})$
 and $\bQ$ for 
 $5^2_1=(8/3)$, $6^2_2=(10/3)$ and $6^2_3=(12/5)$
 appearing in the description of the zeta functions,
 are obtained from $\bQ$ by adjoining
 all the $\bP^1$-coordinates of the degenerate fibers
 of the canonical components
 considered as conic bundles over $\bP^1$.
 We also note that
%
 the trace fields
 (which are same as the invariant trace fields) of
 $5^2_1$, $6^2_2$ and $6^2_3$ are $\bQ(\sqrt{-1})$,
 $\bQ(\sqrt{-3})$ and $\bQ(\sqrt{-7})$ respectively
 (\Cf \cite{GMM2}, \S$5$).
%
 
 Many researches are made on 
 the special values of the Dedekind zeta functions 
 of algebraic number fields at the integer points.
 Especially the values at $s=0,1$ are interesting,
 since those values are described by the intrinsic
 invariants of the number fields
 (class number formula).
 \reftable{SVofZETA} below shows the special values
 (more precisely the coefficient
 of the lowest degree in the Laurent expansion
 when the point is a zero or a pole) of the above
 zeta functions at $s=0,1,2$.
 Note that we have only computed the special values of the
 essential parts of $\zeta(L_i,s)$
 (namely the terms of the products of the Dedekind zeta functions
 up to rational functions in $p^{-s}$).
 In the table
 $R_K$ is the regulator of an algebraic number field $K$.
 Note that
 $R_{\bQ(\sqrt{2})} = \log \left(1 + \sqrt{2}\right)$
 and $R_{\bQ(\sqrt{5})} = \log \left( \dfrac{1 + \sqrt{5}}{2}\right)$ respectively.



\begin{table}[htb]
\begin{center}
\begin{tabular}{c||c|c|c}
\hline
 & $s=0$ & $s=1$ & $s=2$ \\
\hline
 $\zeta(L_0,s)$ & zero, order $1$ & pole, order $1$ & pole, order $3$ \\
 & $\dfrac{-\zeta(3)}{2^{10} 3^3 \pi^2 }$ & $ \dfrac{R_{\bQ(\sqrt{2})}}{2^5 3}$ & $- \dfrac{\sqrt{2}\pi^4 R_{\bQ(\sqrt{2})}}{2^4 3^2}$ \\
\hline
 $\zeta(L_1,s)$ & zero, order $1$ & pole, order $1$ & pole, order $2$ \\
 & $ \dfrac{\zeta(3)}{2^7 3^2 5^2 \pi^2}$ & $- \dfrac{R_{\bQ(\sqrt{5})}^2}{2^4 3}$ & $- \dfrac{ \pi^{6}R_{\bQ(\sqrt{5})}^2}{2^2 3^3 5} $ \\
\hline
 $\zeta(L_2,s)$ & zero, order $1$ & pole, order $2$ & \\
 & $- \dfrac{\zeta(3)}{2^4 \pi^2}$ & $- \dfrac{1}{2^2 3}$ & $- \dfrac{\pi^4}{2^3 3^2}$ \\
\hline
\end{tabular}
\end{center}
\caption{Special values of the main terms of $\zeta(L_i,s)$ at $s=0,1,2$}
\label{SVofZETA}
\end{table}

%
%

\noindent
 It is shown by Smyth (\cite{BoSpe}, appendix or \cite{Smythexp3}, corollary $2$)
 that
 the value $\zeta(3)/\pi^2$ is expressed as
 $2m(1+x+y+z) / 7$,
 where the (logarithmic) Mahler measure $m(P)$ is defined by
$$
 m(P) =
 \int_{0}^1\cdots \int_0^1 \log \left\lv P\left(e^{2\pi t_1\sqrt{-1}},\cdots,e^{2\pi t_n\sqrt{-1}}\right)\right\rv  dt_1\cdots dt_n.
$$
 for a multivariable Laurent polynomial $P \in \bC[x_1^{\pm 1},\cdots,x_n^{\pm 1}]$.
 Note that this expression is not unique,
 namely
 $\zeta(3)/\pi^2$ is also expressed as
 $m\left((1+w)(1+x) + (1-w)(1-x)y\right)/7$ (\cite{Matildeexamples}, \S 2)
 in terms of the Mahler measure of a $4$-variable polynomial.
 Since the other terms are the logarithms of some algebraic integers,
 they are expressed as the Mahler measures of
 those integers.
 Therefore we see that the special values of $\zeta(L_i,s)$
 at $s=0,1,2$ are the products of the Mahler measures
 of certain polynomials up to algebraic numbers and
 some positive powers of $\pi$.

 There is a relation between the volumes of hyperbolic
 $3$-manifolds and the special values of Dedekind zeta
 functions and $L$-functions of imaginary quadratic fields
 (\Cf \cite{BoydMahlerInv,BoRo1}).
 Especially it is known that the volume of the figure $8$
 knot complement in the $3$-sphere is
 equal to the Mahler measure
 of the $A$-polynomial of the figure $8$ knot.
%
%
 Thus it would be interesting to seek for expressions
 of the above special values by the Mahler measures
 of certain link invariants.


 The number of conjugacy classes of $\SL_2$-representations
 of link groups over a fixed finite field is a computable link invariant.
 The explicit formula for the number of conjugacy classes of
 $\SL_2$-representations of torus knot groups over finite fields was given
 by Weiping Li and Liang Xu
 in 2003--2004 (\cite{L-X1,L-X2}).
 Recently the number of $\SL_2(\bFp)$-conjugacy classes
 of all the $\SL_2(\bFp)$-representations has been computed
 for small primes
 by Kitano and Suzuki (\cite{KitSuz}) for all the knots
 in the Rolfsen's knot table,
 which is effective enough to distinguish them.
 The description of the Hasse-Weil zeta functions
 in the Theorem
 provides the exact number of $\SL_2(\bFpn)$-conjugacy classes
 of $\SL_2(\bFpn)$-representations
 (in the canonical components)
 of the arithmetic two-bridge link groups
 for any prime number $p$ and $n\ge 1$.

 In \S $3.1$ of her thesis \cite{EmilyLandesThesis}
 Landes studied
 the infinite family of hyperbolic two-bridge links
 (namely hyperbolic two-bridge links with Schubert normal form $(6n+2/3)$, see also \cite{MarPet} for the hyperbolicity)
 obtained from the magic $3$ manifold by $1/n$ Dehn filling
 on one cusp,
 and showed that for $n=1,2,3,4$
 their canonical components of the character varieties
 have conic bundle structure over $\bP^1$.
 The number of rational points
 and the local zeta functions of standard conic bundles
 over finite fields were determined by
 Tsfasman \cite{Tsfasman} and Rybakov \cite{Rybakov}.
 Hence it would be interesting to
 show whether the canonical components of this family
 have conic bundle structure or not,
 and to study the Hasse-Weil zeta
 functions of this family of hyperbolic two-bridge links.

%
%

 It is well known that
 the Dehn filling on one cusp of the Whitehead link complement
 produces infinitely many hyperbolic $3$ manifolds
 with one cusp.
 In particular we can obtain the twist knot complements
 from the Whitehead link complement by $1/n$ Dehn filling.
 Recently Macasieb, Petersen and van luijk showed
 (\cite{MPL}) that
 the $\SL_2$ character varieties of the twist knots
 are hyperelliptic curves.
 It would be interesting to study
 the zeta functions of twist knots in terms of
 that of the Whitehead link complement,
 which may provide information on
 the Hasse-Weil zeta functions of certain family of
 hyperelliptic curves.

 The author is grateful to the referees for their careful reading
 of the original manuscript.
 Especially he appreciates for pointing out an error in the computation
 of the number of rational points on the complements of degenerate fibers
 of the conic bundles, which fixed the description of the zeta functions
 in Theorem.

\section{Projective models of character varieties of $5^2_1$, $6^2_2$ and $6^2_3$}

 Let $\bP^n:=\bP^n_k$ be
 the $n$-dimensional projective space over
 a field $k$.
 Let
\begin{align*}
 f_0&:=z^3 -xyz^2+(x^2+y^2-2)z -xy,\\
 f_1&:=z^4-xyz^3+(x^2+y^2-3)z^2-xyz+1,\\
 f_2&:=z^3-xyz^2+(x^2+y^2-1)z-xy
\end{align*}
\noindent
 be the defining polynomial of
 the canonical component of the $\SL_2(\bC)$ character variety
 of the Whitehead link complement $5^2_1=(8/3)$,
 $6^2_2=(10/3)$ and $6^2_3=(12/5)$ respectively
 (for details, see \cite{CanATBL}, $1.1$).
 We remark that
 the canonical components of arithmetic two-bridge links
 are hypersurfaces.
 Thus the defining polynomials are uniquely determined
 as monic irreducible polynomials.
 Let 
\begin{align*}
 F:=F_0&:=u^2z^3-xyz^2w +(x^2+y^2-2u^2)zw^2 -xyw^3,\\
 F_1&:=u^2z^4 - xyz^3w + (x^2+y^2-3u^2)z^2w^2 -xyzw^3 +u^2w^4,\\
 F_2&:=u^2z^3-xyz^2w+(x^2+y^2-u^2)zw^2-xyw^3
\end{align*}
 be the corresponding homogeneous polynomials of $f_i$
 in $\bP^2\times \bP^1$, where
$$
 \bP^2\times\bP^1:=\{(x:y:u,\;z:w) \;\mid\;(x:y:u)\in\bP^2,\;(z:w) \in \bP^1\}.
$$
 Let $S_i:=V(F_i) \subset \bP^2\times\bP^1$
 be the projective (singular) surface defined by $F_i$
 over the field $k$.
 For the singularities of $S_i$
 we can show the following lemma by
 the Jacobian criterion
 (we will not use this lemma in what follows in this note).

\begin{lemma}
 Let $k$ be an algebraically closed field with
 characteristic $\Char k=p\ge 0$.
\begin{enumerate}
\item
 When $\Char k \neq 2$,
 the surface $S_0$ has the following $4$ singular points
$$
 (1:0:0,\; 1:0), \; (0:1:0,\; 1:0), \; (1:1:0,\; 1:1),\; (1:-1:0,\; 1:-1).
$$
\noindent
 When $\Char k =2$,
 in addition to the above $3$ points
 (note that $(1:1:0,\; 1:1),\; (1:-1:0,\; 1:-1)$ are same in $\Char k =2$),
 the other singular points of $S_0$ are as follows:
$$
 (x:1:x+1,\; 1:1),\; (x:x+1:1,\; 1:1),\; (1:y:y+1,\; 1:1),\; (0:0:1,\; 0:1),
$$
 where $x$, $y \in k$.

\item
 When $\Char k \neq 2,5$,
 the surface $S_1$ has the following $6$ singular points:
\begin{align*}
 (1:0:0,\;1:0),&\; (0:1:0,\;1:0),\; (1:0:0,\;0:1),\\
 (0:1:0,\;0:1),&\; (1:1:0,\;1:1),\;(1:-1:0,\;1:-1).
\end{align*}
\noindent
 When $\Char k = 5$,
 in addition to the above $6$ points in $\Char k \neq 2,5$ case,
 $S_1$ has $2$ more singular points $(0:0:1,\; 1:\pm 2)$.

\noindent 
 When $\Char k =2$,
 in addition to the $5$ points in $\Char k \neq 2,5$ case,
 $S_1$ has the following singular points:
$$
 (x:x+1:1,\; 1:1),\;(x:1:x+1,\; 1:1),\;
 (1:y:y+1,\; 1:1),\;
 (0:0:1,\; 1:w),
$$
 where $x$, $y \in k$ and
 $w \in k$ satisfies $w^2+w+1=0$.
%
%


\item 
 When $\Char k \neq 2$,
 the projective surface $S_2=V(F_2)$
 has the following $4$ singular points:
$$
 (1:0:0,\;1:0),\; (0:1:0,\;1:0),\; (1:1:0,\;1:1),\; (1:-1:0,\;1:-1).
$$
\noindent
 When $\Char k =2$,
 in addition to the above $3$ points,
 $S_2$ has the following singularities:
$$
 (x:x:1,\; 1:1),\; (1:1:u,\; 1:1)
$$
 where $x,\; u \in k$.

\end{enumerate}
\end{lemma}
 
%
\noindent
 Consider the projection
 $\phi_i:S_i \ra \bP^1$ which is defined by
 $(x:y:u,\;z:w)\mapsto (z:w)$.
 Then the surface $S_i$ can be considered
 as a (singular) conic bundle
 over $\bP^1$, that is,
 each fiber is a (singular) projective curve of degree $2$.


\section{Number of rational points on $S_i$}

 In what follows
 we consider the surfaces $S_i$ over the finite field $\bFp$
 and compute the number of their $k$-rational points,
 where $k=\bFq$ is the finite field
 with $q=p^n$ elements.
 For this purpose
 we regard the set of $k$-rational points of $S_i$
 as the union of fibers
 at the $k$-rational points of $\bP^1$
 and study its degenerate fibers.
 Then we compute the local and global zeta functions of them.


\subsection{$5^2_1$ case}

%
%
%
 When $p\neq 2$
 we see directly from the defining polynomial
 by the Jacobian criterion that
 the surface $S_0$ has $6$ degenerate fibers at 
 $(1:0),\,(0:1),\,(1:\pm 1),\,(1:\pm \tfrac{1}{\sqrt{2}})$.
 In fact, the degenerate fibers of $\phi_0:S_0\ra \bP^1$ are
 expressed as follows:
\begin{align*}
 \phi^{-1}_0(1:0)&=\{(x:y:u,\;1:0)\in \bP^2\times\bP^1 \mid u^2=0\},\\
 \phi^{-1}_0(0:1)&=\{(x:y:u,\;0:1)\in \bP^2\times\bP^1 \mid xy=0\},\\
 \phi^{-1}_0(1:\pm\tfrac{1}{\sqrt{2}})&=\{(x:y:u,\;1:\pm\tfrac{1}{\sqrt{2}})\in \bP^2\times\bP^1 \mid
 \tfrac{1}{2}(x \mp \sqrt{2}y)(x \mp \tfrac{1}{\sqrt{2}}y)=0\},\\
\phi^{-1}_0(1:\pm 1)&=\{(x:y:u,\;1:\pm 1)\in \bP^2\times\bP^1 \mid (x\mp y)-u)((x\mp y)+u)=0\}.
\end{align*}

\noindent
 Then the number of $k=\bFq$-rational points on each fiber
 is as follows:
$$
 \left\lv\phi^{-1}_0(1:0)\right\rv =q + 1,\quad
 \left\lv\phi^{-1}_0(0:1)\right\rv =2q + 1,\quad
 \left\lv\phi^{-1}_0(1:\pm 1)\right\rv =2q + 1,
$$
$$
\left\lv \phi^{-1}_0\left(1:\pm\tfrac{1}{\sqrt{2}}\right)\right\rv =
\begin{cases}
 2q + 1 & \text{ if } \quadres{2}{p} = 1\\
 (2q + 1)\dfrac{1 + (-1)^n}{2} &
 \text{ if } \quadres{2}{p} = -1.
\end{cases}
$$
\noindent
 Here $\tquadres{2}{p}$ is the Legendre symbol of $2$
 modulo $p$.
 Since any other fiber at $(1:w)$
 except these degenerate fibers
 have the $k$-rational point $(1:w:0,\; 1:w)$,
 they are
 isomorphic to $\bP^1$ over $k$
 (\Cf \cite{Liu}, Chap.\;7, Prop.\;4.1).
 Therefore, for $q=p^n$,
 the number $N_{n,0}$ of all the $k=\bFq$-rational points of $S_0$ is
\begin{align*}
 N_{n,0}&=
\begin{cases}
 (q+1)(q-5) \\ (q+1)\left(q-3-(1+(-1)^n)\right)
\end{cases}
 +(q+1) + 3(2q+1) + 2
\begin{cases}
 2q + 1 & \text{ if } \quadres{2}{p} = 1\\
 (2q + 1)\dfrac{1 + (-1)^n}{2} &
 \text{ if } \quadres{2}{p} = -1
\end{cases}\\
    &= 
\begin{cases}
 q^2 + 7q + 1 & \text{ if } \quadres{2}{p} = 1\\
 q^2 + 5q + 1 + q(1 + (-1)^n) &
 \text{ if } \quadres{2}{p} = -1.
\end{cases}
\end{align*}
\noindent
 Therefore, for an arbitrary odd prime $p$,
 the local zeta function $Z_0(p,T)$ is computed as follows:
\begin{align*}
Z_0(p,T)&:=\exp\left(\sum_{n = 1}^\infty
 \dfrac{N_{n,0}}{n} T^n \right)\\
 &=
(1 - p^2T)^{-1}(1 - pT)^{-5}(1 - T)^{-1} \times
\begin{cases}
 (1 - p T)^{-2} & \text{ if } \quadres{2}{p} = 1\\
 \disp{\exp}\left(\sum_{n = 1}^\infty
 \dfrac{2p^{2n}}{2n} T^{2n} \right) & \text{ if } \quadres{2}{p} = -1\end{cases}\\
 &=
\begin{cases}
 (1 - p^2T)^{-1}(1 - pT)^{-7}(1 - T)^{-1} & \text{ if } \quadres{2}{p} = 1\\
 (1 - p^2T)^{-1}(1 - pT)^{-6}(1 - T)^{-1}(1 + pT)^{-1} & \text{ if } \quadres{2}{p} = -1.\end{cases}\\
\end{align*}


%
 When $p=2$,
 the surface $S_0$
 has $3$ degenerate fibers at 
 $(1:0),\,(0:1),\,(1: 1)$.
 In fact the degenerate fibers of $\phi_0:S_0\ra \bP^1$ are
 expressed as follows:
\begin{align*}
 \phi^{-1}_0(1:0)&=\{(x:y:u,\;1:0)\in \bP^2\times\bP^1 \mid u^2=0\},\\
 \phi^{-1}_0(0:1)&=\{(x:y:u,\;0:1)\in \bP^2\times\bP^1 \mid xy=0\},\\
\phi^{-1}_0(1: 1)&=\{(x:y:u,\;1:1)\in \bP^2\times\bP^1 \mid (x + y + u)^2=0\}.
\end{align*}
\noindent
 Hence we can compute the number of $k$-rational points
 on each degenerate fiber.
 The result is
$$
 \left\lv \phi^{-1}_0(1:0) \right\rv = q + 1, \quad
 \left\lv \phi^{-1}_0(0:1) \right\rv = 2q + 1, \quad
 \left\lv \phi^{-1}_0(1: 1) \right\rv = q + 1.
$$
\noindent
 Since any other fiber except these degenerate fibers is
 isomorphic to $\bP^1$ over $k$,
 the number $N_{n,0}$ of all the $k=\bFq$-rational
 points of $S_0$ for $p=2$ is
\begin{align*}
 N_{n,0}&= (q+1)(q-2) +2(q+1) + (2q+1)\\
    &= q^2 + 3q + 1.
\end{align*}
 Thus the local zeta function $Z_0(2,T)$ is
\begin{align*}
Z_0(2,T)&=\exp\left(\sum_{n = 1}^\infty
 \dfrac{N_{n,0}}{n} T^n \right)
 = (1 - 2^2T)^{-1}(1 - 2T)^{-3}(1 - T)^{-1}.
\end{align*}
%


 Thus the Hasse-Weil zeta function
$$
 \zeta(S_0,s):= \prod_{p} Z_0(p, p^{-s})
$$
 is written as follows:
\begin{align*}
\zeta(S_0,s)&=\left(\prod_p (1 - p^{2-s})^{-1}(1 - p^{1-s})^{-6}(1 - p^{-s})^{-1} \right) (1 - 2^{1-s})^{3} \times\\
& \prod_{\tquadres{2}{p}=1}(1 - p^{1-s})^{-1}
\times \prod_{\tquadres{2}{p}=-1}(1 + p^{1-s})^{-1}\\
&= \zeta(s-2)\zeta(s-1)^6 \zeta(s) L\left(\tquadres{2}{\cdot},s-1\right)(1 - 2^{1-s})^{3}.
\end{align*}

\noindent 
 Here 
$$
 \dquadres{2}{\cdot} : (\bZ/8\bZ)^{\times} \ra \bC^{\times}; \quad n \mapsto \dquadres{2}{n}
$$
 is the Dirichlet character
 associated with the Legendre symbol $\dquadres{2}{\cdot}$
 and
$$
 L\left(\tquadres{2}{\cdot},s\right)=\sum_{n\ge 1}\dfrac{\tquadres{2}{n}}{n^s}=\prod_p \left(1 - \tquadres{2}{p}p^{-s}\right)^{-1}
$$
 is the Dirichlet $L$-function of $\tquadres{2}{\cdot}$.


\subsubsection{non-affine part}
 To have the description of the Hasse-Weil zeta function
 $\zeta(f_0,s)$
 here we compute the number of $k$-rational points of
 the non-affine part $V(F_0)\ssm V(f_0)$
 and its zeta function.
 We see that
 the set $V(F_0)\ssm V(f_0)$ consists of points which have
 coordinates $(x:y:0,\;1:0)$ or $(x:y:0,\;z:1)$.
 Indeed it consists of three subsets
\begin{align*}
 V(F_0)\ssm V(f_0)(k)&=\{(x:y:0,\;1:0)\mid (x:y)\in \bP^1_k \}\\
 &\cup\;\{(z:w:0,\;z:w)\mid (z:w)\in\bP^1_k \}\;\cup\;\{(w:z:0,\;z:w)\mid (z:w)\in \bP^1_k \}.
\end{align*}
\noindent
 Then it is easy to compute the number $N'_{n,0}$ of $k$-rational points
 of $V(F_0)\ssm V(f_0)$ and its zeta functions.
$$
 N'_{n,0} =
 \begin{cases}
 3q - 1, & \text{ if } p \neq 2,\\
 3q    , & \text{ if } p = 2.
 \end{cases}
$$
$$
 Z'_0(p,T) = \exp\left(\sum_{n = 1}^\infty
 \dfrac{N'_{n,0}}{n} T^n \right)=
 \begin{cases}
 (1 -pT)^{-3}(1 - T),& \text{ if } p \neq 2,\\
 (1 - 2T)^{-3},& \text{ if } p =2.
\end{cases}
$$
$$
 \tilde{\zeta}(S_0,s):=\prod_p Z'_0(p,p^{-s})= \zeta(s-1)^3 \zeta(s)^{-1}(1-2^{-s})^{-1}.
$$
 Thus we have the description of the Hasse-Weil zeta function
 of the canonical component $V(f_0)$ of $5^2_1$.
\begin{align*}
 \zeta(5^2_1,s):=\zeta(f_0,s)&:=\prod_p \exp\left(\sum_{n=1}^{\infty}\dfrac{\#V(f_0)(\bFpn)}{n} (p^{-s})^n\right)= \zeta(S_0,s)/\tilde{\zeta}(S_0,s)\\
             &=\zeta(s-2)\zeta(s-1)^3 \zeta(s)^2 L\left(\tquadres{2}{\cdot},s-1\right)(1 - 2^{1-s})^3 (1 - 2^{-s})\\
             &=\zeta(s-2)\zeta(s-1)^2 \zeta(s)^2 \zeta_{\bQ(\sqrt{2})}(s-1)(1 - 2^{1-s})^3 (1 - 2^{-s}).
\end{align*}



\subsection{$6^2_3$ case}


 Since $6^2_3$ case is similar to the previous $5^2_1$ case
 and the computation is rather simple,
 we treat this case next and summarize the result without
 much detail.
 When $p\neq 2$
 the surface $S_2$ has $4$ degenerate fibers at 
 $(1:0),\,(0:1),\,(1:\pm 1)$.
 In fact, the degenerate fibers of $\phi_2:S_2\ra \bP^1$ are
 expressed as follows:
\begin{align*}
 \phi^{-1}_2(1:0)&=\{(x:y:u,\;1:0)\in \bP^2\times\bP^1 \mid u^2=0\}=\{(x:y:0,\;1:0)\},\\
 \phi^{-1}_2(0:1)&=\{(x:y:u,\;0:1)\in \bP^2\times\bP^1 \mid xy=0\}=\{(x:0:u,\;0:1)\}\vee\{(0:y:u,\;0:1)\},\\
\phi^{-1}_2(1:\pm 1)&=\{(x:y:u,\;1:\pm 1)\in \bP^2\times\bP^1 \mid (x\mp y)^2=0\}=\{(x:\pm x: u,\;1:\pm 1)\in \bP^2\times\bP^1\}.
\end{align*}
\noindent
 Hence the number of $k=\bFq$-rational points on each fiber
 is as follows:
\begin{align*}
 \left\lv \phi^{-1}_2(1:0) \right\rv &=q + 1,\\
 \left\lv \phi^{-1}_2(0:1) \right\rv &=2q + 1,\\
 \left\lv \phi^{-1}_2(1:\pm 1) \right\rv &=q + 1.
\end{align*}
\noindent
 Since any other fiber except these degenerate fibers are
 isomorphic to $\bP^1$ over $k$,
 the number $N_{n,2}$ of all the $k$-rational points of $S_2$ is
\begin{align*}
 N_{n,2}&=
 (q+1)(q-3) +(q+1) + (2q+1) + 2(q+1)\\
    &= q^2+3q+1.
\end{align*}


 For the surface $S_2$
 we can see that the number of $k$-rational points
 for $p=2$ is
 same as the odd prime case.
 In fact, when $p=2$,
 the surface $S_2$
 has $3$ degenerate fibers at 
 $(1:0),\,(0:1),\,(1: 1)$.
 The fibers at these points are same as the ones in
 the $p\neq 2$ case.
 Hence we can compute the number of $k$-rational points
 on each degenerate fiber.
 The result is
$$
 \left\lv \phi^{-1}_2(1:0) \right\rv = q + 1, \quad
 \left\lv \phi^{-1}_2(0:1) \right\rv = 2q + 1, \quad
 \left\lv \phi^{-1}_2(1: 1) \right\rv = q + 1.
$$
\noindent
 Since any other fiber except these degenerate fibers is
 isomorphic to $\bP^1$ over $k$,
 the number $N_{n,2}$ of all the $k$-rational
 points of $S_2$ for $p=2$ is
\begin{align*}
 N_{n,2}&= (q+1)(q-2) +2(q+1) + (2q+1)\\
    &= q^2 + 3q + 1.
\end{align*}

\noindent
 Therefore, for an arbitrary prime $p$ (including $2$),
 the local zeta function $Z_2(p,T)$ is computed as follows:
\begin{align*}
Z_2(p,T)&:=\exp\left(\sum_{n = 1}^\infty
 \dfrac{N_{n,2}}{n} T^n \right)= (1 - p^2T)^{-1}(1 - pT)^{-3}(1 - T)^{-1}.
\end{align*}



 Thus the Hasse-Weil zeta function
$$
 \zeta(S_2,s):= \prod_{p} Z_2(p, p^{-s})
$$
 is written as follows:
\begin{align*}
\zeta(S_2,s)&=\prod_p (1 - p^{2-s})^{-1}(1 - p^{1-s})^{-3}(1 - p^{-s})^{-1}\\
&= \zeta(s-2)\zeta(s-1)^3 \zeta(s).
\end{align*}

\subsubsection{non-affine part}
 The set $V(F_2)\ssm V(f_2)$ consists of points which have
 coordinates $(x:y:0,\;1:0)$ or $(x:y:0,\;z:1)$.
 Indeed it consists of three subsets
\begin{align*}
 V(F_2)\ssm V(f_2)(k)&=\{(x:y:0,\;1:0)\mid (x:y)\in \bP^1_k \}\\
 &\cup\;\{(z:w:0,\;z:w)\mid (z:w)\in\bP^1_k \}\;\cup\;\{(w:z:0,\;z:w)\mid (z:w)\in \bP^1_k \}.
\end{align*}
\noindent
 Then it is easy to compute the number $N'_{n,2}$ of $k$-rational points
 of $V(F_2)\ssm V(f_2)$ and its zeta functions.
$$
 N'_{n,2} =
 \begin{cases}
 3q - 1, & \text{ if } p \neq 2,\\
 3q    , & \text{ if } p = 2.
 \end{cases}
$$
$$
 Z'_2(p,T) = \exp\left(\sum_{n = 1}^\infty
 \dfrac{N'_{n,2}}{n} T^n \right)=
 \begin{cases}
 (1 -pT)^{-3}(1 - T),& \text{ if } p \neq 2,\\
 (1 - 2T)^{-3},& \text{ if } p = 2.
\end{cases}
$$
$$
 \tilde{\zeta}(S_2,s):=\prod_p Z'_2(p,p^{-s})= \zeta(s-1)^3 \zeta(s)^{-1}(1-2^{-s})^{-1}.
$$
 Thus we have the description of the Hasse-Weil zeta function
 of the canonical component $V(f_2)$ of $6^2_3$.
\begin{align*}
 \zeta(6^2_3,s):=\zeta(f_2,s)&:=\prod_p \exp\left(\sum_{n=1}^{\infty}\dfrac{\#V(f_2)(\bFpn)}{n} (p^{-s})^n\right)= \zeta(S_2,s)/\tilde{\zeta}(S_2,s)\\
             &= \zeta(s-2)\zeta(s)^2 (1-2^{-s}).
\end{align*}


\subsection{$6^2_2$ case}

%


 In this case
 we first assume $\Char k \neq 2,5$.
 Then the surface $S_1$
 has $8$ degenerate fibers at 
 $(1:0),\,(0:1),\,(1:\pm 1),\,(1:\pm \tfrac{\sqrt{5}\pm 1}{2})$.

\begin{align*}
 \phi^{-1}_1(1:0)&=\{(x:y:u,\;1:0)\in \bP^2\times\bP^1 \mid u^2=0\},\\
 \phi^{-1}_1(0:1)&=\{(x:y:u,\;0:1)\in \bP^2\times\bP^1 \mid u^2=0\},\\
 \phi^{-1}_1\left(1:\tfrac{\sqrt{5}\pm 1}{2}\right)&=\{(x:y:u,\;1:\tfrac{\sqrt{5}\pm 1}{2})\in \bP^2\times\bP^1 \mid
 \tfrac{3\pm\sqrt{5}}{2}(x - \tfrac{\sqrt{5}+1}{2}y)(x - \tfrac{\sqrt{5}-1}{2}y)=0\},\\
 \phi^{-1}_1 \left(1:-\tfrac{\sqrt{5}\pm 1}{2}\right)&=\{(x:y:u,\;1:-\tfrac{\sqrt{5}\pm 1}{2})\in \bP^2\times\bP^1 \mid
 \tfrac{3\pm\sqrt{5}}{2}(x + \tfrac{\sqrt{5}+1}{2}y)(x + \tfrac{\sqrt{5}-1}{2}y)=0\},\\
\phi^{-1}_1(1:\pm 1)&=\{(x:y:u,\;1:\pm 1)\in \bP^2\times\bP^1 \mid ((x\mp y)-u)((x\mp y)+u)=0\}.
\end{align*}
\noindent
 Hence the number of $k=\bFq$-rational points on each fiber
 for $p \neq 2,5$
 is as follows:
%
%
$$
 \left\lv \phi^{-1}_1(1:0) \right\rv =q + 1,\quad
 \left\lv \phi^{-1}_1(0:1) \right\rv =q + 1,\quad
 \left\lv \phi^{-1}_1(1:\pm 1) \right\rv =2q + 1,
$$
\begin{align*}
 \left\lv \phi^{-1}_1\left(1:\pm \tfrac{\sqrt{5}\pm 1}{2}\right) \right\rv &=
\begin{cases}
 2q + 1 & \text{ if } \quadres{5}{p} = 1\\
 (2q + 1)\dfrac{1 + (-1)^n}{2} &
 \text{ if } \quadres{5}{p} = -1.
\end{cases}
\end{align*}
\noindent
 Here $\tquadres{5}{p}$ is the Legendre symbol of $5$
 modulo $p$.
 Since any other fiber except these degenerate fibers are
 isomorphic to $\bP^1$ over $k$,
 the number $N_{n,1}$ of all the $k$-rational points of $S_1$ is
\begin{align*}
 N_{n,1}&=
\begin{cases} (q+1)(q-7) \\ (q+1)\left(q-3-2(1+(-1)^n)\right)
\end{cases}
 +2(q+1) + 2(2q+1) + 4
\begin{cases}
 2q + 1 & \text{ if } \quadres{5}{p} = 1\\
 (2q + 1)\dfrac{1 + (-1)^n}{2} &
 \text{ if } \quadres{5}{p} = -1
\end{cases}\\
    &= 
\begin{cases}
 q^2 + 8q + 1 & \text{ if } \quadres{5}{p} = 1,\\
 q^2 + 4q + 1 + 2q(1 + (-1)^n) &
 \text{ if } \quadres{5}{p} = -1.
\end{cases}
\end{align*}
\noindent
 Therefore, for an arbitrary prime $p \neq 2,5$,
 the local zeta function $Z_1(p,T)$ is computed as follows:
\begin{align*}
Z_1(p,T)&:=\exp\left(\sum_{n = 1}^\infty
 \dfrac{N_{n,1}}{n} T^n \right)\\
 &= (1 - p^2T)^{-1}(1 - pT)^{-4}(1 - T)^{-1} \times
\begin{cases}
 (1 - p T)^{-4} & \text{ if } \quadres{5}{p} = 1\\
 \disp{\exp}\left(\sum_{n = 1}^\infty
 \dfrac{4p^{2n}}{2n} T^{2n} \right) & \text{ if } \quadres{5}{p} = -1 \end{cases}\\
 &= 
\begin{cases}
 (1 - p^2T)^{-1}(1 - pT)^{-8}(1 - T)^{-1} & \text{ if } \quadres{5}{p} = 1\\
 (1 - p^2T)^{-1}(1 - pT)^{-4}(1 - T)^{-1}(1 - p^2T^2)^{-2} & \text{ if } \quadres{5}{p} = -1.\end{cases}\\
\end{align*}


 When $p=2$,
 the surface $S_1$
 has $5$ degenerate fibers at 
 $(1:0),\,(0:1),\,(1: 1), (1:\ga^{\pm 1})$,
 where $\ga$ is a root of $T^2+T+1=0$.
 In fact, the degenerate fibers of $\phi_1:S_1\ra \bP^1$ are
 expressed as follows:
\begin{align*}
 \phi^{-1}_1(1:0)&=\{(x:y:u,\;1:0)\in \bP^2\times\bP^1 \mid u^2=0\},\\
 \phi^{-1}_1(0:1)&=\{(x:y:u,\;0:1)\in \bP^2\times\bP^1 \mid xy=0\},\\
\phi^{-1}_1(1: 1)&=\{(x:y:u,\;1:1)\in \bP^2\times\bP^1 \mid (x + y + u)^2=0\},\\
\phi^{-1}_1(1: \ga^{\pm 1})&=\{(x:y:u,\;1:\ga^{\pm 1})\in \bP^2\times\bP^1 \mid (x - \ga y)(x - \ga^{-1}y)=0\}.
\end{align*}
\noindent
 Hence we can compute the number of $k$-rational points
 on each degenerate fiber.
 The result is
$$
 \left\lv \phi^{-1}_1(1:0) \right\rv = q + 1, \quad
 \left\lv \phi^{-1}_1(0:1) \right\rv = q + 1, \quad
$$
$$
 \left\lv \phi^{-1}_1(1: 1) \right\rv = q + 1, \quad
 \left\lv \phi^{-1}_1(1: \ga^{\pm 1}) \right\rv = (2q + 1)\dfrac{1 + (-1)^n}{2}.
$$
\noindent
 Since any other fiber except these degenerate fibers is
 isomorphic to $\bP^1$ over $k$,
 the number $N_{n,1}$ of all the $k$-rational
 points of $S_1$ for $p=2$ is
\begin{align*}
 N_{n,1}&= (q+1)\left(q-2 - (1+(-1)^n)\right) +3(q+1) + 2(2q+1)\dfrac{1 + (-1)^n}{2}\\
    &= q^2 +2q + 1 + q\left(1 + (-1)^n \right).
\end{align*}
 Thus the local zeta function $Z_1(2,T)$ is
\begin{align*}
Z_1(2,T)&=\exp\left(\sum_{n = 1}^\infty
 \dfrac{N_{n,1}}{n} T^n \right)
 = (1 - 2^2T)^{-1}(1 - 2 T)^{-2}(1 - T)^{-1}(1 - 2^2T^2)^{-1}.
\end{align*}


 When $p=5$,
 the surface $S_1$
 has $6$ degenerate fibers at 
 $(1:0),\,(0:1),\,(1: \pm 1), (1: \pm 2)$.
 In fact, the degenerate fibers of $\phi_1:S_1\ra \bP^1$ are
 expressed as follows:
\begin{align*}
 \phi^{-1}_1(1:0)&=\{(x:y:u,\;1:0)\in \bP^2\times\bP^1 \mid u^2=0\},\\
 \phi^{-1}_1(0:1)&=\{(x:y:u,\;0:1)\in \bP^2\times\bP^1 \mid xy=0\},\\
\phi^{-1}_1(1:\pm 1)&=\{(x:y:u,\;1:\pm 1)\in \bP^2\times\bP^1 \mid ((x\mp y)-u)((x\mp y)+u)=0\},\\
\phi^{-1}_1(1: \pm 2)&=\{(x:y:u,\;1:\pm 2)\in \bP^2\times\bP^1 \mid (x - 2 y)(x + 2y)=0\}.
\end{align*}
\noindent
 Hence we can compute the number of $k$-rational points
 on each degenerate fiber.
 The result is
$$
 \left\lv \phi^{-1}_1(1:0) \right\rv = q + 1, \quad
 \left\lv \phi^{-1}_1(0:1) \right\rv = q + 1, \quad
 \left\lv \phi^{-1}_1(1: \pm 1) \right\rv = 2q + 1, \quad
 \left\lv \phi^{-1}_1(1: \pm 2) \right\rv = 2q + 1.
$$
\noindent
 Since any other fiber except these degenerate fibers is
 isomorphic to $\bP^1$,
 the number $N_{n,1}$ of all the $k$-rational
 points of $S_1$ for $p=5$ is
\begin{align*}
 N_{n,1}&= (q+1)(q-5) + 2(q+1) + 4(2q+1)= q^2 + 6q + 1.
\end{align*}
 Thus the local zeta function $Z_1(5,T)$ is
\begin{align*}
Z_1(5,T)&=\exp\left(\sum_{n = 1}^\infty
 \dfrac{N_{n,1}}{n} T^n \right)
 = (1 - 5^2T)^{-1}(1 - 5T)^{-6}(1 - T)^{-1}.
\end{align*}
%

\noindent
 Hence the Hasse-Weil zeta function
$$
 \zeta(S_1,s):= \prod_{p} Z_1(p, p^{-s})
$$
 is written as follows:
\begin{align*}
\zeta(S_1,s)&=
\zeta(s-2)\zeta(s-1)^6\zeta(s)
 L\left(\tquadres{5}{\cdot},s-1\right)^2\times 
 (1 - 2^{1-s})^{3}(1 + 2^{1-s}).
\end{align*}

\subsubsection{non-affine part}
%
 The set $V(F_1)\ssm V(f_1)$ consists of four subsets
\begin{align*}
 V(F_1)\ssm V(f_1)(k)&=\{(x:y:0,\;1:0)\mid (x:y)\in \bP^1_k \}\;\cup\;\{(x:y:0,\;0:1)\mid (x:y)\in \bP^1_k \}\\
 &\cup\;\{(z:w:0,\;z:w)\mid (z:w)\in\bP^1_k \}\;\cup\;\{(w:z:0,\;z:w)\mid (z:w)\in \bP^1_k \}.
\end{align*}
\noindent
 Then it is easy to compute the number $N'_{n,1}$ of $k$-rational points
 of $V(F_1)\ssm V(f_1)$ and its zeta functions.
\begin{align*}
 N'_{n,1} &= 2(q+1)+ (q-1) +
 \begin{cases}
 q - 3, & \text{ if } p \neq 2\\
 q - 2    , & \text{ if } p = 2
 \end{cases}\\
 &=4q -
 \begin{cases}
  2, & \text{ if } p \neq 2,\\
  1    , & \text{ if } p = 2.
 \end{cases}
\end{align*}
$$
 Z'_1(p,T) =
 \begin{cases}
 (1 -pT)^{-4}(1 - T)^2,& \text{ if } p \neq 2,\\
 (1 - 2T)^{-4}(1 - T),& \text{ if } p =2.
\end{cases}
$$
$$
 \tilde{\zeta}(S_1,s):=\prod_p Z'_1(p,p^{-s})= \zeta(s-1)^4 \zeta(s)^{-2}(1-2^{-s})^{-1}.
$$
 Thus we have the description of the Hasse-Weil zeta function
 of the canonical component $V(f_1)$ of $6^2_2$.
\begin{align*}
 \zeta(6^2_2,s):=\zeta(f_1,s)&:=\prod_p \exp\left(\sum_{n=1}^{\infty}\dfrac{\#V(f_1)(\bFpn)}{n} (p^{-s})^n\right)= \zeta(S_1,s)/\tilde{\zeta}(S_1,s)\\
             &=
 \zeta(s-2)\zeta(s-1)^2\zeta(s)^3  L\left(\tquadres{5}{\cdot},s-1\right)^2
 \times 
(1 - 2^{1-s})^{3}(1 + 2^{1-s})(1 - 2^{-s})
\\
%
             &=
 \zeta_{\bQ(\sqrt{5})}(s-1)^2 \zeta(s-2)\zeta(s)^3 
 \times 
(1 - 2^{1-s})^{3}(1 + 2^{1-s})(1 - 2^{-s}).
\end{align*}



\begin{thebibliography}{10}

\bibitem{BoSpe}
D.~W.\ Boyd, \emph{Speculations concerning the range of {M}ahler's measure},
  Canad.\ Math.\ Bull.\ \textbf{24} (1981), 453--469.


\bibitem{BoydMahlerInv}
\bysame, \emph{Mahler's measure and invariants of hyperbolic manifolds}, Number
  theory for the millennium, {I} ({U}rbana, {IL}, 2000),
  MA, 2002, 127--143.


\bibitem{BoRo1}
D.~W.\ Boyd and F.\ Rodriguez-Villegas, \emph{Mahler's measure and the
  dilogarithm. {I}}, Canad.\ J.\ Math.\ \textbf{54} (2002), 468--492.


\bibitem{GMM2}
F.~W.\ Gehring, C.~Maclachlan, and G.~J.\ Martin, \emph{Two-generator arithmetic
  {K}leinian groups. {II}}, Bull.\ London Math.\ Soc.\ \textbf{30} (1998), 258--266.


\bibitem{CanATBL}
S.\ Harada, \emph{Canonical components of character varieties of arithmetic
  two-bridge link complements}, preprint, http://arxiv.org/abs/1112.3441.

\bibitem{HWfigHarada}
\bysame, \emph{Hasse-{W}eil zeta function of absolutely irreducible {${\rm
  SL}_2$}-representations of the figure 8 knot group}, Proc.\ Amer.\ Math.\ Soc.\ \textbf{139} (2011), 3115--3125.


\bibitem{KitSuz}
T.\ Kitano and M.\ Suzuki, \emph{On the number of {$SL(2;\Bbb Z/p\Bbb
  Z)$}-representations of knot groups}, J.\ Knot Theory Ramifications
  \textbf{21} (2012), 1250003, 18.

\bibitem{Matildeexamples}
M.~N.\ Lal{\'{\i}}n, \emph{Some examples of {M}ahler measures as multiple
  polylogarithms}, J.\ Number Theory \textbf{103} (2003), 85--108.


\bibitem{EmilyLandesThesis}
E.\ Landes, \emph{On the canonical components of character varieties of
  hyperbolic 2-bridge link complements}, Available at
  http://repositories.lib.utexas.edu/bitstream/handle/2152/ETD-UT-2011-08-2877/LANDES-DISSERTATION.pdf?sequence=1.

\bibitem{EmilyLandes}
\bysame, \emph{Identifying the canonical component for the {W}hitehead link},
  Math.\ Res.\ Lett.\ \textbf{18} (2011), 715--731.


\bibitem{L-X1}
W.\ Li and L.\ Xu, \emph{Counting {${\rm SL}\sb 2({\bf F}\sb {2\sp s})$}
  representations of torus knot groups}, Acta Math.\ Sin.\ (Engl. Ser.) \textbf{19} (2003), 233--244.


\bibitem{L-X2}
\bysame, \emph{Counting {${\rm SL}\sb 2({\bf F}\sb q)$}-representations of
  torus knot groups}, J.\ Knot Theory Ramifications \textbf{13} (2004), 401--426.


\bibitem{Liu}
Q.\ Liu, \emph{Algebraic geometry and arithmetic curves}, Oxford Grad.\ Texts in Math., \textbf{6}, Oxford Univ.\ Press, Oxford, 2002.


\bibitem{MPL}
M.~L.\ Macasieb, K.~L.\ Petersen, and R.~M.\ van Luijk, \emph{On
  character varieties of two-bridge knot groups}, Proc.\ Lond.\ Math.\ Soc.\ (3)
  \textbf{103} (2011), 473--507.


\bibitem{BookMaRe}
C.\ Maclachlan and A.~W.\ Reid, \emph{The arithmetic of hyperbolic
  3-manifolds}, Grad.\ Texts in Math., \textbf{219}, Springer-Verlag, New
  York, 2003.


\bibitem{MarPet}
B.\ Martelli and C.\ Petronio, \emph{Dehn filling of the ``magic''
  3-manifold}, Comm.\ Anal.\ Geom.\ \textbf{14} (2006), 969--1026.


\bibitem{Rybakov}
S.\ Rybakov, \emph{Zeta functions of conic bundles and {D}el {P}ezzo
  surfaces of degree 4 over finite fields}, Mosc.\ Math.\ J.\ \textbf{5} (2005), 919--926, 974.


\bibitem{Smythexp3}
C.~J.\ Smyth, \emph{An explicit formula for the {M}ahler measure of a family
  of 3-variable polynomials}, J.\ Th\'eor.\ Nombres Bordeaux \textbf{14} (2002), 683--700.


\bibitem{Tsfasman}
M.~A.\ Tsfasman, \emph{Nombre de points des surfaces sur un corps fini},
  Arithmetic, geometry and coding theory ({L}uminy, 1993), de Gruyter, Berlin,
  1996, 209--224.

\end{thebibliography}

\providecommand{\bysame}{\leavevmode\hbox to3em{\hrulefill}\thinspace}
\providecommand{\MR}{\relax\ifhmode\unskip\space\fi MR }
\providecommand{\MRhref}[2]{%
  \href{http://www.ams.org/mathscinet-getitem?mr=#1}{#2}
}
\providecommand{\href}[2]{#2}


\vspace{0.5cm}

%
%

\noindent
Previous address:\\
Shinya Harada\\
School of Mathematics\\
Korea Institute for Advanced Study (KIAS)\\
Hoegiro 85, Dongdaemun-gu\\
 Seoul 130-722, KOREA\\
{\tt harada@kias.re.kr}

\vspace{0.5cm}

\noindent
Current address:\\
Shinya Harada\\
Graduate School of Information Science and Engineering,\\
 Tokyo Institute of Technology,\\
W8-37 2-12-1 Ookayama, Meguro-ku, Tokyo 152-8550, Japan\\
{\tt harada.s.al@m.titech.ac.jp}

\end{document}